\documentclass[12pt]{article}
\usepackage{amsthm,amsfonts}
\usepackage{fullpage}

\newtheorem{theorem}{Theorem}

\title{Words avoiding reversed subwords}
\author{Narad Rampersad and Jeffrey Shallit\\
School of Computer Science \\
University of Waterloo \\
Waterloo, ON, N2L 3G1 \\
CANADA \\
{\tt nrampersad@math.uwaterloo.ca} \\
{\tt shallit@graceland.math.uwaterloo.ca}}

\begin{document}
\date{\today}
\maketitle

\begin{abstract}
We examine words $w$ satisfying the following property:  if $x$ is a
subword of $w$ and $|x|$ is at least $k$ for some fixed $k$, then the
reversal of $x$ is not a subword of $w$.
\end{abstract}

\section{Introduction}
Let $\Sigma$ be a finite, nonempty set called an \emph{alphabet}.
We denote the set of all finite words over the alphabet
$\Sigma$ by $\Sigma^*$.  The empty word is represented by $\epsilon$.
Let $\Sigma_k$ denote the alphabet $\{0,1,\ldots,k-1\}$.

Let $\mathbb{N}$ denote the set $\{0,1,2,\ldots\}$.
An \emph{infinite word} is a map from $\mathbb{N}$ to $\Sigma$.
The set of all infinite words over the alphabet $\Sigma$ is
denoted $\Sigma^\omega$.

A map $h:\Sigma^* \rightarrow \Delta^*$ is called a \emph{morphism} if
$h$ satisfies $h(xy)=h(x)h(y)$ for all $x,y\in\Sigma^*$.  A morphism
may be defined by specifying its action on $\Sigma$.  Morphisms
may also be applied to infinite words in the natural way.

If $w \in \Sigma^*$ is written $w = w_1 w_2 \cdots w_n$, where each
$w_i \in \Sigma$, then the \emph{reversal} of $w$, denoted $w^R$,
is the word $w_n w_{n-1} \cdots w_1$.

If $y$ is a nonempty word, then the word $yyy\cdots$ is written as
$y^\omega$.  If an infinite word $\mathbf{w}$ can be written in the form
$y^\omega$ for some nonempty $y$, then $\mathbf{w}$ is said to be
\emph{periodic}.  If $\mathbf{w}$ can be written in the form $y'y^\omega$
for some nonempty $y$, then $\mathbf{w}$ is said to be
\emph{ultimately periodic}.

A \emph{square} is a word of the form $xx$, where $x\in\Sigma^*$ is nonempty.
A word $w'$ is called a \emph{subword} (resp. a \emph{prefix} or a
\emph{suffix}) of $w$ if $w$ can be written in the form $uw'v$
(resp. $w'v$ or $uw'$) for some $u,v\in\Sigma^*$.
We say a word $w$ is \emph{squarefree} (or \emph{avoids squares})
if no subword of $w$ is a square.

\section{Avoiding reversed subwords}
Szilard \cite{Szi03} has asked the following question:
\begin{quotation}
Does there exist an infinite word $\mathbf{w}$ such that if $x$ is a
subword of $\mathbf{w}$, then $x^R$ is not a subword of $\mathbf{w}$?
\end{quotation}

Clearly there must be some restriction on the length of $x$: if $|x|=1$,
then all nonempty words fail to have the desired property.  For
$|x| \geq 2$, however, we have the following result.

\begin{theorem}
\label{tern_per}
There exists an infinite word $\mathbf{w}$ over $\Sigma_3$ such that if
$x$ is a subword of $\mathbf{w}$ and $|x| \geq 2$, then $x^R$ is not a
subword of $\mathbf{w}$.  Furthermore, $\mathbf{w}$ is unique up to
permutation of the alphabet symbols.
\end{theorem}

\begin{proof}
Note that if $|x| \geq 3$ and both $x$ and $x^R$ are subwords of
$\mathbf{w}$, then there is a prefix $x'$ of $x$ such that $|x'|=2$
and $(x')^R$ is a suffix of $x^R$.  Hence it suffices to show
the theorem for $|x|=2$.  We show that the infinite word
$$\mathbf{w} = (012)^\omega = 012012012012 \cdots$$
has the desired property.  To see this, consider the set $\mathcal{A}$
consisting of all subwords of $\mathbf{w}$ of length two.  We have
$\mathcal{A} = \{01,12,20\}$.
Noting that if $x \in \mathcal{A}$, then $x^R \not\in \mathcal{A}$,
we conclude that if $x$ is a subword of $\mathbf{w}$ and $|x| \geq 2$,
then $x^R$ is not a subword of $\mathbf{w}$.

To see that $\mathbf{w}$ is unique up to permutation of the alphabet
symbols, consider another word $\mathbf{w}'$ satisfying the conditions
of the theorem, and suppose that $\mathbf{w}'$ begins with 01.  Then
01 must be followed by 2, 12 must be followed by 0, and 20 must be
followed by 1.  Hence,
$$\mathbf{w}' = (012)^\omega = 012012012012 \cdots = \mathbf{w}.$$
\end{proof}

Note that the solution given in the proof of Theorem~\ref{tern_per} is
periodic.  In the following theorem, we give a nonperiodic solution to
this problem for $|x| \geq 3$.

\begin{theorem}
There exists an infinite nonperiodic word $\mathbf{w}$ over
$\Sigma_3$ such that if $x$ is a subword of $\mathbf{w}$ and $|x| \geq 3$,
then $x^R$ is not a subword of $\mathbf{w}$.
\end{theorem}

\begin{proof}
By reasoning similar to that given in the proof of Theorem~\ref{tern_per},
it suffices to show the theorem for $|x|=3$.
Let $\mathbf{w'}$ be an infinite nonperiodic word over $\Sigma_2$.  For
example, if $\mathbf{w'} = 11010010001\cdots$, then $\mathbf{w'}$
is nonperiodic. Define the morphism
$h : \Sigma_2^\omega \rightarrow \Sigma_3^\omega$ by
\begin{eqnarray*}
0 & \rightarrow & 0012 \\
1 & \rightarrow & 0112.
\end{eqnarray*}
Then $\mathbf{w} = h(\mathbf{w'})$ has the desired property.
Consider the set $\mathcal{A}$ consisting of all subwords of
$\mathbf{w}$ of length three.  We have
$$\mathcal{A} = \{001, 011, 012, 112, 120, 200, 201\}.$$
Noting that if $x \in \mathcal{A}$, then $x^R \not\in \mathcal{A}$,
we conclude that if $x$ is a subword of $\mathbf{w}$ and $|x| \geq 3$,
then $x^R$ is not a subword of $\mathbf{w}$.

To see that $\mathbf{w}$ is not periodic, suppose the contrary; \emph{i.e.},
suppose that $\mathbf{w} = y^\omega$ for some $y \in \Sigma_3^*$.  Clearly,
$|y| > 4$.  Suppose then that $y$ begins with $h(0)$.  Noting that the only
way to obtain 00 from $h(ab)$, where $a,b \in \Sigma_2$, is as a prefix of
$h(0)$, we see that $y = h(y')$ for some $y' \in \Sigma_2^*$.  Hence,
$\mathbf{w} = \left(h\left(y'\right)\right)^\omega =
h\left(\left(y'\right)^\omega\right)$, and so $\mathbf{w'} = (y')^\omega$
is periodic, contrary to our choice of $\mathbf{w'}$.
\end{proof}

Over a two letter alphabet we have the following negative result.

\begin{theorem}
Let $k \leq 4$ and let $w$ be a word over $\Sigma_2$ such that
if $x$ is a subword of $w$ and $|x| \geq k$, then $x^R$ is not
a subword of $w$.  Then $|w| \leq 8$.
\end{theorem}

\begin{proof}
As mentioned previously, if $k=1$ the result holds trivially.  If $k=2$,
note that all binary words of length at least three must contain one of the
following words: 00, 11, 010, or 101.  Similarly, if $k=3$, note that all
binary words of length at least five must contain one of the following words:
000, 010, 101, 111, 0110, or 1001; and if $k=4$, note that all binary words of
length at least nine must contain one of the following words: 0000, 0110,
1001, 1111, 00100, 01010, 01110, 10001, 10101, or 11011.  Hence, $|w| \leq 8$,
as required.
\end{proof}

For $|x| \geq 5$, however, we find that there \emph{are} infinite words with
the desired property.

\begin{theorem}
\label{geq5}
There exists an infinite word $\mathbf{w}$ over $\Sigma_2$ such that if
$x$ is a subword of $\mathbf{w}$ and $|x| \geq 5$, then $x^R$ is not a
subword of $\mathbf{w}$.
\end{theorem}

\begin{proof}
By reasoning similar to that given in the proof of Theorem~\ref{tern_per},
it suffices to show the theorem for $|x|=5$.  We show that the infinite word
$$\mathbf{w} = (001011)^\omega = 001011001011001011 \cdots$$
has the desired property.  To see this, consider the set $\mathcal{A}$
consisting of all subwords of $\mathbf{w}$ of length five.  We have
$$\mathcal{A} = \{00101, 01011, 01100, 10010, 10110, 11001\}.$$
Noting that if $x \in \mathcal{A}$, then $x^R \not\in \mathcal{A}$,
we conclude that if $x$ is a subword of $\mathbf{w}$ and $|x| \geq 5$,
then $x^R$ is not a subword of $\mathbf{w}$.
\end{proof}

Let $z$ be the word 001011.  We denote the \emph{complement}
of $z$ by $\bar{z}$, \emph{i.e.}, the word obtained by
substituting 0 for 1 and 1 for 0 in $z$.  Let $\mathcal{B}$
be the set defined as follows:
$$\mathcal{B} = \{x \,|\, x \mbox{ is a cyclic shift of }\, z \mbox{ or }
\bar{z}\}.$$  We have the following characterization of the words
satisfying the conditions of Theorem~\ref{geq5}.

\begin{theorem}
\label{per}
Let $\mathbf{w}$ be an infinite word over $\Sigma_2$ such that if
$x$ is a subword of $\mathbf{w}$ and $|x| \geq 5$, then $x^R$ is not a
subword of $\mathbf{w}$.  Then $\mathbf{w}$ is ultimately periodic.
Specifically, $\mathbf{w}$ is of the form $y'y^\omega$, where
$y' \in \{\epsilon,0,1,00,11\}$ and $y \in \mathcal{B}$.
\end{theorem}

\begin{proof}
By reasoning similar to that given in the proof of Theorem~\ref{tern_per},
it suffices to show the theorem for $|x|=5$.  We call a word $w \in \Sigma_2^*$
\emph{valid} if $w$ satisfies the property that if
$x$ is a subword of $w$ and $|x|=5$, then $x^R$ is not a
subword of $w$.  We have the following two facts,
which may be verified computationally.
\begin{enumerate}
\item All valid words of length 9 are of the form $y'yy''$, where
$y' \in \{\epsilon,0,1,00,11\}$, $y \in \mathcal{B}$, and $y'' \in \Sigma_2^*$.
\item Let $w$ be a valid word of the form $yy''$, where
$y \in \mathcal{B}$ and $y'' \in \Sigma_2^*$.  Then if $|w|=15$,
$y$ is a prefix of $y''$.
\end{enumerate}

We will prove by induction on $n$ that for all $n \geq 1$,
$y'y^n$ is a prefix of $\mathbf{w}$, where $y' \in \{\epsilon,0,1,00,11\}$
and $y \in \mathcal{B}$.

If $n=1$, then by applying the first fact to the prefix of $\mathbf{w}$
of length 9, we have that $y'y$ is a prefix of $\mathbf{w}$, as required.

Assume then that $y'y^n$ is a prefix of $\mathbf{w}$.  We can thus write
$\mathbf{w} = y'y^{n-1}y\mathbf{w'}$, for some
$\mathbf{w'} \in \Sigma_2^\omega$.  By applying the second fact to
the prefix of $y\mathbf{w'}$ of length 15, we have that $y$ is a
prefix of $\mathbf{w'}$.  Hence $\mathbf{w} = y'y^{n-1}yy\mathbf{w''} =
y'y^{n+1}\mathbf{w''}$, for some $\mathbf{w''} \in \Sigma_2^\omega$,
as required.

We therefore conclude that if $\mathbf{w}$ satisfies the
conditions of the theorem, then $\mathbf{w}$ is of the form $y'y^\omega$,
where $y' \in \{\epsilon,0,1,00,11\}$ and $y \in \mathcal{B}$.
\end{proof}

Next we give a nonperiodic solution to this problem for $|x| \geq 6$.

\begin{theorem}
There exists an infinite nonperiodic word $\mathbf{w}$ over
$\Sigma_2$ such that if $x$ is a subword of $\mathbf{w}$ and $|x| \geq 6$,
then $x^R$ is not a subword of $\mathbf{w}$.
\end{theorem}

\begin{proof}
By reasoning similar to that given in the proof of Theorem~\ref{tern_per},
it suffices to show the theorem for $|x|=6$.
Let $\mathbf{w'}$ be an infinite nonperiodic word over $\Sigma_2$.
Define the morphism $h : \Sigma_2^\omega \rightarrow \Sigma_2^\omega$ by
\begin{eqnarray*}
0 & \rightarrow & 0001011 \\
1 & \rightarrow & 0010111.
\end{eqnarray*}
We show that the infinite word $\mathbf{w} = h(\mathbf{w'})$
has the desired property.  To see this, consider the set $\mathcal{A}$
consisting of all subwords of $\mathbf{w}$ of length six.  We have
\begin{eqnarray*}
\mathcal{A} & = & \{000101, 001011, 010110, 010111, 011000, 011001, 011100, \\
& & 100010, 100101, 101100, 101110, 110001, 110010, 111000, 111001\}.
\end{eqnarray*}
Noting that if $x \in \mathcal{A}$, then $x^R \not\in \mathcal{A}$,
we conclude that if $x$ is a subword of $\mathbf{w}$ and $|x| \geq 6$,
then $x^R$ is not a subword of $\mathbf{w}$.

To see that $\mathbf{w}$ is not periodic, suppose the contrary; \emph{i.e.},
suppose that $\mathbf{w} = y^\omega$ for some $y \in \Sigma_2^*$.  Clearly,
$|y| > 7$.  Suppose then that $y$ begins with $h(0)$.  Noting that the only
way to obtain 000 from $h(ab)$, where $a,b \in \Sigma_2$, is as a prefix of
$h(0)$, we see that $y = h(y')$ for some $y' \in \Sigma_2^*$.  Hence,
$\mathbf{w} = \left(h\left(y'\right)\right)^\omega =
h\left(\left(y'\right)^\omega\right)$, and so $\mathbf{w'} = (y')^\omega$
is periodic, contrary to our choice of $\mathbf{w'}$.
\end{proof}

Finally we consider words avoiding squares as well as reversed subwords.
It is easy to check that no binary word of length $\geq 4$ avoids
squares.  However, Thue \cite{Thu06} gave an example of a infinite
squarefree ternary word.  Over a four letter alphabet we have the
following negative result, which may be verified computationally.

\begin{theorem}
\label{sqf4}
Let $w$ be a squarefree word over $\Sigma_4$ such that
if $x$ is a subword of $w$ and $|x| \geq 2$, then $x^R$ is not
a subword of $w$.  Then $|w| \leq 20$.
\end{theorem}

In contrast with the result of Theorem~\ref{sqf4},
Alon \emph{et al.} \cite{AGHR02} have noted that over a four letter alphabet
there exists an infinite squarefree word that avoids palindromes $x$ where
$|x| \geq 2$.  (A \emph{palindrome} is a word $x$ such that $x = x^R$.)
However, over a five letter alphabet there are infinite words with an
even stronger avoidance property.

\begin{theorem}
There exists an infinite squarefree word $\mathbf{w}$ over $\Sigma_5$
such that if $x$ is a subword of $\mathbf{w}$ and $|x| \geq 2$, then
$x^R$ is not a subword of $\mathbf{w}$.
\end{theorem}

\begin{proof}
By reasoning similar to that given in the proof of Theorem~\ref{tern_per},
it suffices to show the theorem for $|x|=2$.
Let $\mathbf{w'}$ be an infinite squarefree word over $\Sigma_3$.
Define the morphism $h : \Sigma_3^\omega \rightarrow \Sigma_5^\omega$ by
\begin{eqnarray*}
0 & \rightarrow & 012 \\
1 & \rightarrow & 013 \\
2 & \rightarrow & 014.
\end{eqnarray*}
We show that the infinite word $\mathbf{w} = h(\mathbf{w'})$
has the desired property.

First we note that to verify that $\mathbf{w}$ is squarefree, it suffices
by a theorem of Thue \cite{Thu12} (see also \cite{BEM79}, \cite{Ber79}, and
\cite{Cr82}) to verify that $h(w)$ is squarefree for all 12 squarefree words
$w \in \Sigma_3^*$ such that $|w|=3$.  This is left to the reader.

To see that if $x$ is a subword of $\mathbf{w}$ and $|x|=2$, then
$x^R$ is not a subword of $\mathbf{w}$, consider the set $\mathcal{A}$
consisting of all subwords of $\mathbf{w}$ of length 2.  We have
$$\mathcal{A} = \{01, 12, 13, 14, 20, 30, 40\}.$$
Noting that if $x \in \mathcal{A}$, then $x^R \not\in \mathcal{A}$,
we conclude that if $x$ is a subword of $\mathbf{w}$ and $|x| \geq 2$,
then $x^R$ is not a subword of $\mathbf{w}$.
\end{proof}

\end{document}